\documentclass[12pt, reqno]{amsart}
\usepackage{amsfonts}
\usepackage{amsmath}
\usepackage{amssymb}
\usepackage{eufrak}
\usepackage{mathrsfs}
\usepackage{bbding}

\numberwithin{equation}{section} 

\input xy
\xyoption{all}
\title{On the existence of certain weak Fano threefolds of Picard number two}

\author{Maxim Arap} \author{Joseph Cutrone} \author{Nicholas Marshburn}

\begin{document}

\maketitle

\newcommand{\bP}{\mathbb{P}}
\newcommand{\bC}{\mathbb{C}}
\newcommand{\mc}{\mathcal}
\newcommand{\ra}{\rightarrow}
\newcommand{\thra}{\twoheadrightarrow}
\newcommand{\mb}{\mathbb}
\newcommand{\mrm}{\mathrm}
\newcommand{\p}{\prime}
\newcommand{\ms}{\mathscr}
\newcommand{\pl}{\partial}
\newcommand{\ti}{\tilde}
\newcommand{\wti}{\widetilde}
\newcommand{\ol}{\overline}
\newcommand{\ul}{\underline}
\newcommand{\Sp}{\mrm{Spec}\,}

\renewcommand{\thefootnote}{}

\newtheorem{pro1}{Proposition}[section]
\newtheorem{pro2}[pro1]{Proposition}
\newtheorem{lem0.5}[pro1]{Lemma}
\newtheorem{lem0}[pro1]{Lemma}
\newtheorem{pro3}[pro1]{Proposition}
\newtheorem{pro4}[pro1]{Proposition}
\newtheorem{pro5}[pro1]{Proposition}
\newtheorem{pro6}[pro1]{Proposition}
\newtheorem{pro7}[pro1]{Proposition}
\newtheorem{pro8}[pro1]{Proposition}

\newtheorem{def1}[pro1]{Definition}

\newtheorem{thm0.1}[pro1]{Theorem}
\newtheorem{thm0.2}[pro1]{Theorem}
\newtheorem{thm0.3}[pro1]{Theorem}
\newtheorem{thm0.4}[pro1]{Theorem}
\newtheorem{thm0.5}[pro1]{Theorem}

\newtheorem{thm1}[pro1]{Theorem}
\newtheorem{thm2}[pro1]{Theorem}
\newtheorem{thm3}[pro1]{Theorem}
\newtheorem{thm4}[pro1]{Theorem}
\newtheorem{thm5}[pro1]{Theorem}

\newtheorem{cor1}[pro1]{Corollary}
\newtheorem{cor2}[pro1]{Corollary}
\newtheorem{cor3}[pro1]{Corollary}
\newtheorem{cor5}[pro1]{Corollary}
\newtheorem{cor6}[pro1]{Corollary}

\newtheorem{lem1}[pro1]{Lemma}
\newtheorem{lem2}[pro1]{Lemma}
\newtheorem{lem3}[pro1]{Lemma}
\newtheorem{Prokhorov's lemma}[pro1]{Lemma}
\newtheorem{lem4}[pro1]{Lemma}
\newtheorem{lem5}[pro1]{Lemma}
\newtheorem{lem6}[pro1]{Lemma}
\newtheorem{lem7}[pro1]{Lemma}
\newtheorem{lem8}[pro1]{Lemma}
\newtheorem{lem9}[pro1]{Lemma}
\newtheorem{lem10}[pro1]{Lemma}
\newtheorem{lem11}[pro1]{Lemma}

\newtheorem{rem0}[pro1]{Remark}
\newtheorem{rem1}[pro1]{Remark}
\newtheorem{rem2}[pro1]{Remark}
\newtheorem{rem2.5}[pro1]{Remark}
\newtheorem{rem3}[pro1]{Remark}
\newtheorem{rem4}[pro1]{Remark}
\newtheorem{rem5}[pro1]{Remark}

\footnote{\Envelope\, Maxim Arap: marap@math.jhu.edu; Joseph Cutrone: jcutrone@towson.edu; Nicholas Marshburn: marshbur@math.jhu.edu.}

\begin{abstract} 
This article settles the question of existence of smooth weak Fano threefolds of Picard number two with small anti-canonical map and previously classified numerical invariants obtained by blowing up certain curves on smooth Fano threefolds of Picard number 1 with the exception of 12 numerical cases. 
\end{abstract} 

\thispagestyle{empty}
\section{Introduction}
A Fano variety is a smooth projective variety whose anti-canonical class is ample. To study birational maps between Fano varieties it is interesting to relax the positivity condition on the anti-canonical class and require it to be nef and big but not ample.  Projective varieties whose anti-canonical class is nef and big but not ample are called weak Fano or almost Fano. In what follows, all varieties are defined over $\mb{C}$ and all Fano and weak Fano varieties are assumed to be smooth, unless stated otherwise.

Fano threefolds have been classified mainly through the works of Fano, Iskovskikh, Shokurov, Mori and Mukai (\cite{Isk77}, \cite{Isk78}, \cite{Sh79}, \cite{MM84}). As a next step, it is interesting to classify birational maps between Fano threefolds. This question fits into the general framework of Sarkisov program. Within this framework, weak Fano threefolds of Picard number two yield birational maps called Sarkisov links between Fano threefolds of Picard number one. This is one of the motivations to classify weak Fano threefolds of Picard number two.

Besides the contribution to the classification of Sarkisov links, the results of this paper have applications to the construction of certain $G_2$-manifolds and Calabi-Yau threefolds, see \cite{CHNP12} and \cite{CHNP13}.  

In \cite{JPR05}, \cite{Tak09}, \cite{JPR11} and \cite{CM13} one finds the numerical constraints that weak Fano threefolds of Picard number two must satisfy, obtaining a finite list of possibilities, some of which were shown to exist while others were left as numerical possibilities. The recent article \cite{BL12} constructs some of the Sarkisov links which were previously not known to exist. The aim of this article is to study further the geometric realizability of the numerical possibilities which were left open in \cite{CM13}.

Let $X$ be a weak Fano threefold of Picard number two such that the anti-canonical system $|{-}K_X|$ is free and gives a small contraction $\psi\colon X \ra X^\p$. By \cite{Ko89}, the $K_X$-trivial curves can be flopped. More precisely, there is a commutative diagram
$$\xymatrix{
X \ar@{-->}[rr]^-\chi \ar[d]_-{\phi} \ar[rd]^-\psi & & X^+ \ar[ld]_-{\psi^+} \ar[d]^-{\phi^+} \\
Y & X^\p & Y^+,
}$$
where $\chi$ is an isomorphism outside of the exceptional locus of $\psi$. The numerical cases, which we study here, are of E1-E1 type, i.e., both $\phi$ and $\phi^+$ are assumed to be divisorial contractions of type E1 in the sense of \cite{Mo82}. In particular, $Y$ is a smooth Fano variety of Picard number one and $\phi$ is the blow-up of $Y$ along a smooth irreducible curve $C$, whose degree and genus are denoted by $d$ and $g$, respectively. Likewise $\phi^+$ is the blow-up of $Y^+$ along a smooth irreducible curve $C^+$ of degree $d^+$ and genus $g^+$.  The exceptional divisors of the blow-ups $\phi$ and $\phi^+$ are denoted by $E$ and $E^+$, respectively.

\begin{def1}\label{def1}
For a smooth Fano threefold $Y$ of Picard number 1, define $\mathscr{C}(Y)$ to be the set of all $(d,g)$ such that there is a numerical link of type E1-E1 as above with $\phi$ being the blow-up of a smooth irreducible curve $C \subset Y$ of degree $d$ and genus $g$.  
\end{def1}

The numerical classification of Sarkisov links of E1-E1 type with $|{-}K_X|$ giving a small contraction can be found in \cite[Table E1-E1]{CM13}. The numerical cases of type E1-E1, whose geometric realizability was left open in \cite{CM13} are listed at the end of this paper for convenience. The main results of this article are as follows. 
 
\begin{thm0.2}[=Theorem \ref{thm2}]
The numerical invariants listed in cases 44-46, 48, 70, 71, 86, 88, 97, 105 in \cite[Table E1-E1]{CM13} of Sarkisov links starting with a smooth quadric are geometrically realizable.  
\end{thm0.2}

\begin{thm0.3}[=Theorem \ref{thm3}]
The numerical invariants listed in cases 30, 34, 37, 41, 64, 84 in \cite[Table E1-E1]{CM13} of Sarkisov links starting with a del Pezzo threefold of degree 4 are geometrically realizable.  
\end{thm0.3}

\begin{thm0.1}[=Theorem \ref{thm1}]
The numerical invariants listed in cases 31, 35, 38, 40, 42, 65, 68, 69, 81, 83, 94, 96, 101 in \cite[Table E1-E1]{CM13} of Sarkisov links starting with a del Pezzo threefold of degree 5 are geometrically realizable. The numerical case 43 is not realizable. 
\end{thm0.1}

\begin{thm0.4}[=Theorem \ref{thm4}]
The numerical invariants listed in cases 3-8, 11-15, 18-20, 23, 55-57, 79, 93, 100, 104 in \cite[Table E1-E1]{CM13} of Sarkisov links starting with a Fano threefold of index 1 are geometrically realizable. The numerical cases 16, 17, 21, 22, 24-26, 58, 60 are not realizable. 
\end{thm0.4}

To summarize, this article settles geometric realizability of all open numerical cases from \cite[Table E1-E1]{CM13} with the exception of 12 cases (see Remarks \ref{rem1}, \ref{rem2}, \ref{rem4} for the explanation of why these cases are left open). The results of the above theorems are summarized in Tables \ref{tb:E1E1Table1}-\ref{tb:E1E1Table3} at the end of this article. 

Throughout the article we shall use the following notation. The index of $Y$ is the largest integer $r$ such that there exists an ample divisor $H$ with ${-}K_Y = rH$. Throughout the article we shall refer to $H$ as the hyperplane class on $Y$. Abusing notation for the sake of simplicity, the pull-back of $H$ to $X$ is also denoted by $H$. The anti-canonical class of $X$ is ${-}K_X = rH-E$.  We let $l$ and $f$ in $N_1(X)$ denote the classes of the pull-back of a general line on $Y$ and a $\phi$-exceptional curve, respectively. The elements of $N_1(X)$ can be written as $ml-nf$ for some $n,m\in \mb{Z}$ and the Mori cone $\ol{\mrm{NE}}(X)$ has two extremal rays $f$ and $r$. The extremal ray $r$ has slope 
$$\mrm{sup}\big\{n/m\,:\,  ml-nf \text{ can be represented by an effective curve}\big\}.$$
If $ml-nf$ is represented by a curve $\ti{\Gamma}$, then $\ti{\Gamma}$ is the proper transform of a curve $\Gamma\subset Y$ such that $\deg \Gamma = m$ and $C\cdot \Gamma =n$. The intersection pairing is given by the formulas
$$H\cdot l = 1, E\cdot f=-1, H\cdot f = E\cdot l = 0.$$

\section{Preliminaries} \label{Prelim}
The following proposition, suggested to us by the referee as a simplification of the exposition, will be used throughout the paper to check freeness of the anti-canonical system $|{-}K_X|$. 

\begin{pro2} \label{pro2}
Let $S$ be a smooth K3 with $\mrm{Pic}(S) =\mathbb{Z}H \oplus \mathbb{Z}C$, with $H$ very
ample and $C$ a smooth (irreducible) curve. Assume $H^2=2n$, $C \cdot H=d$ and
$C^2=2(g-1)$. Let $k>0$ be an integer. Then $kH-C$ is nef if and only if

\vspace{0.2cm}

(*) $2nk>d, nk^2-dk+g-1 \geq 0$ and $(2nk-d,nk^2-dk+g) \neq (2n+1,n+1)$.

\vspace{0.2cm}

Furthermore, $kH-C$ is free if and only if it is nef and we
are not in the case

\vspace{0.2cm}

(**) $d^2-4n(g-1)=1$, and $2nk-d-1$ or $2nk-d+1$ divides $2n$.

\vspace{0.2cm}

\end{pro2}

\begin{proof} Set $C_k:=kH-C, d_k:=2nk-d$ and $g_k:=nk^2-dk+g$. Then $C_k \cdot H=d_k$,
$(C_k)^2=2(g_k-1)$ and the two first conditions in (*) are equivalent to
$C_k \cdot H>0$ and $C_k^2 \geq 0$, which are necessary conditions for $C_k$ to be
nef, and imply that $C_k>0$. Then we may conclude that $|C_k|$ is free by
\cite[Prop. 4.4]{Knu02} (with $C_k$ in the role of $C$), and the result about
nefness by the first part of the proof of \cite[Prop. 4.4]{Knu02}. (The case $\Gamma \cdot H=0$ does not
occur as $H$ is assumed to be ample in our proposition.)
\end{proof}

\begin{rem0} In the proof of [Knu02, Prop. 4.4] the divisors contradicting
nefness or global generation of $kH-C$ are explicitly constructed as
linear combinations of $H$ and $C_k$. Hence, on any K3 with hyperplane
section $H$ and a smooth curve $C$ (without any assumption on its Picard
group) with $n,d,g$, as above, $kH-C$ is not nef if $(2nk-d,nk^2-dk+g)
=(2n+1,n+1)$ and not globally generated if (**) holds.
\end{rem0}

\section{Blow-ups of smooth quadric threefolds} \label{Q}
Consider the list of pairs 
$$\ms{C}(Q) = \{(9,2), (10,5), (11,8), (12,11), (13,14), (9,3), (8,0), (10,6), (8,1), (9,4), (8,2), (8,3), (7,1)\},$$
which correspond to the numerical possibilities of $(d,g)$ in cases 44-48, 70-72, 86, 88, 97, 102, 105 in \cite{CM13}. The purpose of this section is to prove that 10 of the above 13 numerical possibilities are geometrically realizable. 

Throughout this section $C_{d,g}\subset \bP^4$ denotes a smooth non-degenerate curve of degree $d$ and genus $g$. For simplicity, we write $C$ instead of $C_{d,g}$, when $d$ and $g$ are understood from the context. A smooth complete intersection of a quadric and a cubic in $\bP^4$ is denoted by $S$. The restriction of $H$ to $S$ is denoted by $H_S$. 

\begin{pro3} \label{pro3}
For each $(d,g) \in \ms{C}(Q)$, there exists a smooth curve $C$ of degree $d$ and genus $g$ on a smooth quadric threefold $Y$ such that the blow-up of $Y$ along $C$ is a smooth weak Fano threefold of Picard number two. 
\end{pro3}
\begin{proof} By \cite[Th.1.1, p.202]{Knu02}, for each $(d,g) \in \ms{C}(Q)$ there exists a smooth curve $C$ of degree $d$ and genus $g$ lying on a smooth K3 surface $S \subset \bP^4$ of degree $6$ with the property $\mrm{Pic}(S) = \mb{Z}H_S\oplus \mb{Z}C$. 

Let $I_S$ be the ideal sheaf of $S$ in $\bP^4$. We have $h^0(\mc{O}_S(2))=14$ and $h^0(\mc{O}_{\bP^4}(2))=15$. Therefore, from the long exact cohomology sequence associated to $0 \ra I_S(2) \ra \mc{O}_{\bP^4}(2) \ra \mc{O}_S(2)\ra 0$ we obtain $h^0(I_S(2))\ge 1$. Since $\deg S=6$ this implies that $S$ is contained in a unique quadric $Y$, \cite[Ex.VIII.14, p.97]{Bea96}. If $Y$ is singular, then $Y$ contains a plane, which cuts out a cubic plane curve $\Gamma$ on $S$ with $\Gamma^2=0$ and $\Gamma\cdot H=3$. We may check that for each $(d,g) \in \ms{C}(Q)$, the class of $\Gamma$ cannot be expressed as an integral linear combination of $H$ and $C$. Therefore, the quadric $Y$ must be smooth. 

Let $\ti{S}\subset X$ be the birational transform of $S$ under the blow-up $\phi\colon X\ra Y$ along $C$. We may and shall identify $\ti{S}$ with $S$. Since $\ti{S} \in |{-}K_X|=|3H-E|$, then $\mrm{Bs}(|{-}K_X|) \subset \mrm{Bs}(|3H-E|_{|\ti{S}})=\mrm{Bs}(|3H_S-C|)$. Using Proposition \ref{pro2}, we may check that the divisor $3H_S-C$ is free on $S$, and therefore, ${-}K_X$ is also free. In particular, ${-}K_X$ is nef. 

To see that ${-}K_X$ is not ample, it suffices to show that $C$ has a trisecant line. Indeed, any trisecant $L$ to $C$ is necessarily contained in the quadric $Y$ and the proper transform of $L$ is $({-}K_X)$-trivial. A formula of Berzolari, see for example \cite{LeB82}, gives that the number of trisecant lines to a curve of degree $d$ and genus $g$ in $\bP^4$ is
$$\theta(d,g)=\binom{d-2}{3}-g(d-4).$$   
We may check that for each $(d,g) \in \ms{C}(Q)$, the number $\theta(d,g)$ is positive. Therefore, $C$ has a trisecant line.

Since ${-}K_X$ is nef, then to show that ${-}K_X$ is big, it suffices to check that $({-}K_X)^3>0$, \cite[Thm.2.2.16, p.144]{LazI}. For this, we use the formula 
$$({-}K_X)^3=({-}K_Y)^3+2K_Y\cdot C -2+2g = 52-6d+2g,$$
see for example \cite[Lem.2.4]{BL12}, and check that for each $(d,g) \in \ms{C}(Q)$ we have $({-}K_X)^3>0$. 
\end{proof}

\begin{thm2} \label{thm2}
The numerical invariants listed in cases 44-46, 48, 70, 71, 86, 88, 97, 105 in \cite[Table E1-E1]{CM13} of Sarkisov links are geometrically realizable.  
\end{thm2}
\begin{proof}
In each of the listed cases the weak Fano threefold $X$ exists by Proposition \ref{pro3}. It suffices to check that $|{-}K_X|$ determines a small contraction and $\phi^+$ is of type E1. The numerical possibilities with $|{-}K_X|$ giving a divisorial contraction are classified in \cite[Table A.4, p.629]{JPR05}. If the morphism $\phi^+$ were not of type E1, then the listed numerical possibilities would have appeared either in \cite[7.4, 7.7, p.486]{JPR11} or in the non-E1-E1 tables in \cite{CM13}.  
\end{proof}

\begin{rem1}\label{rem1} 
The links with numerical invariants listed under 47, 72, 102 in \cite[Table E1-E1]{CM13} also appear in the table \cite[A4, p.629]{JPR05}. It appears to be a delicate question whether $|{-}K_X|$ determines a small contraction in the examples constructed above. The authors are currently investigating these cases and hope to resolve them in a future paper. 
\end{rem1}

\section{Blow-ups of the intersection of two quadrics in $\bP^5$} \label{V4}
Consider the list $$\ms{C}(V_4) =\{(7,0), (8,2), (9,4), (10,6), (11,8), (7,1), (8,3), (7,2)\}$$
of pairs $(d,g)$, which correspond to numerical possibilities 30, 34, 37, 39, 41, 64, 67, 84 in \cite{CM13}. In this section $Y\subset \bP^5$ is a smooth intersection of two quadrics and $S \subset \bP^5$ is a smooth complete intersection of three quadrics. The remaining notation is as in Section \ref{Q}. 

\begin{pro4}
For each $(d,g) \in \ms{C}(V_4)$, there exists a smooth curve $C$ of degree $d$ and genus $g$ on a smooth intersection $Y$ of two quadrics in $\bP^5$ such that the blow-up of $Y$ along $C$ is a smooth weak Fano threefold of Picard number two. 
\end{pro4}
\begin{proof} By \cite[Th.1.1, p.202]{Knu02}, for each $(d,g) \in \ms{C}(V_4)$ there exists a smooth curve $C$ of degree $d$ and genus $g$ lying on a smooth complete intersection K3 surface $S \subset \bP^5$ of degree $8$ with the property $\mrm{Pic}(S) = \mb{Z}H_S\oplus \mb{Z}C$. 

The linear system $|2H-S|$ of quadrics containing $S$ has dimension two and its base locus is $S$. Since $S$ is non-singular, no quadric in $|2H-S|$ can be singular at a point of $S$. By Bertini's theorem, a general member of $|2H-S|$ is also smooth outside of $S$. Therefore, a general member $Q \in |2H-S|$ is a smooth quadric. The restricted linear system $|2H-S|_{|Q}$ is a pencil on $Q$ with base locus $S$. Again by Bertini's theorem, for a general quadric $Q^\p \in |2H-S|$ the intersection $Q \cap Q^\p$ is smooth outside of $S$. Since $S$ is a complete intersection of quadrics, then $Q \cap Q^\p$ must be smooth along $S$ as well. This shows that $S$ is contained in a smooth intersection of two quadrics in $\bP^5$, which is the promised threefold $Y$. 

Let $X\ra Y$ the blow-up  of $C$ and let $\ti{S}$ be the birational transform of $S$. Using Proposition \ref{pro2}, we may check that the divisor $2H_S-C$ is free on $S$. Since $\ti{S} \in |{-}K_X|=|2H-E|$, then, as in the proof of Proposition \ref{pro3}, we conclude that ${-}K_X$ is also free and, in particular, nef. Since ${-}K_X$ is nef and $({-}K_X)^3>0$, then ${-}K_X$ is big. By classification \cite[12.3, p.217]{IP99},  ${-}K_X$ is not ample.
\end{proof}

\begin{thm3} \label{thm3}
The numerical invariants listed in cases 30, 34, 37, 41, 64, 84 in \cite[Table E1-E1]{CM13} of Sarkisov links are geometrically realizable.  
\end{thm3}
\begin{proof}
The proof is the same as that of Theorem \ref{thm2}. 
\end{proof}

\begin{rem2}\label{rem2} 
The links with numerical invariants listed under 39 and 67 in \cite[Table E1-E1]{CM13} also appear in the table \cite[A4, p.629]{JPR05}. It appears to be a delicate question whether $|{-}K_X|$ determines a small contraction in the examples constructed above. The authors are currently investigating these cases and hope to resolve them in a future paper. 
\end{rem2}

\section{Blow-ups of $V_5$} \label{V5}
In this section $Y\subset \bP^9$ is a smooth section of the Pl\"ucker-embedded Grassmannian $\mb{G}(1,4) \subset \bP^9$ of lines in $\bP^4$ by a linear subspace of codimension three. The remaining notation is as before. The open cases are 31, 35, 38, 40, 42, 43, 65, 68, 69, 81, 83, 94, 96, 101 with $(d,g)$ in 
\begin{align}\ms{C}(V_5)=\{(9,0), (10,2), (11,4), (12,6), (13,8), (14,10), (9,1), (10,3)&, (12,7),(9,2), (8,0), \nonumber \\ &(9,3), (8,1), (7,0)\} \nonumber.
\end{align}
In this section we shall show that all of the above numerical cases are geometrically realizable, except the case with $(d,g)=(14,10)$, which is not realizable. 

In the sequel, unless otherwise stated, we fix $(d,g)\in \ms{C}(V_5) \backslash\{(12,7),(13,8),(14,10)\}$. Let $S_{d,g}\subset \bP^6$ be a smooth K3 surface of degree 10 and genus $6$ with the following properties. The surface $S_{d,g}$ is cut out by quadrics and $\mrm{Pic}(S_{d,g}) = \mb{Z}T \oplus \mb{Z} C_{d,g}$, where $T$ is a general hyperplane section of $S_{d,g}$ (a smooth canonical curve of genus $6$) and $C_{d,g}$ is a smooth curve of degree $d$ and genus $g$. The existence of such surfaces was established in \cite{Knu02}.  When there is no danger of confusion we shall omit the subscript $d,g$ from $S_{d,g}$ and $C_{d,g}$ for simplicity.

We shall use Gushel's method (see \cite{Gu82} and \cite{Gu92}) to show that $S$ can be embedded into a smooth del Pezzo threefold $V_5$. This method relies on Maruyama's construction of regular vector bundles, see \cite{Ma73}.  First, we shall construct a globally generated rank two vector bundle $M$ on $S$ such that $h^0(M)=5$. Let $\bP^1_S := \bP^1\times S$ be the trivial $\bP^1$-bundle over $S$ and let $D$ be a $g^1_4$ on $T$. Since $S$ is cut out by quadrics, $T$ has Clifford index $>1$ by \cite{S-D74}, hence $T$ is not trigonal and not a plane quintic. In particular, this implies that $|D|$ is free. 
Fix a point $p \in \bP^1$ and let $Y$ be a general member of the linear system $|\{p\}\times T + \bP^1\times D|$ on $\bP^1_T$. The situation is summarized in the following diagram
$$\xymatrix{
Y \ar@{^{(}->}[r] \ar[rd] & \bP^1_T \ar@{^{(}->}[r] \ar[d] & \bP^1_S \ar[d]^-\pi \\
& T \ar@{^{(}->}[r] & S,
}$$
where the square is Cartesian and $\pi$ is the natural projection. Let $I_Y$ be the ideal sheaf of $Y$ in $\bP^1_S$ and let $G := \{p\}\times S +\bP^1\times T$. Define $M = \pi_*(I_Y\otimes \mc{O}_{\bP^1_S}(G))$. By \cite[Principle 2.6, p.112]{Ma73} or \cite[Lem.1.4]{Gu82}, $M$ is a vector bundle of rank two on $S$. 

Let us show that $M$ is globally generated. Let $W = H^0(D)$. 
As in \cite[1.5.2]{Gu82}, there is a short exact sequence 
\begin{equation} \label{MG}
0 \ra W^*\otimes \mc{O}_S \ra M \ra \mc{O}_T(K_T-D) \ra 0 
\end{equation} 
(note that $T^2$ in Gushel's notation is a hyperplane section of $T$, which in our case is a canonical divisor $K_T$). Using the facts $h^1(W^*\otimes \mc{O}_S) = h^1(\mc{O}_S) = 0$, $\dim W^* = 2$, and $h^0(K_T-D)=3$, 
the long exact sequence  associated  to (\ref{MG}) gives $h^0(M)=5$. To show that $M$ is globally generated, it suffices to prove that $|K_T-D|$ is free. The linear system $|K_T-D|$ is a $g^2_6$. Since $T$ is not a plane quintic, $|K_T-D|$ is free and this
completes the proof that $M$ is globally generated. 

Therefore, the evaluation homomorphism $H^0(M)\otimes \mc{O}_S \twoheadrightarrow M$ determines a morphism $\alpha\colon S \ra \mb{G}(1,4)$. By \cite[Cor.2.19.1, p.121]{Ma73}, $c_1(M)=T$ and $c_2(M)=D$.  Since $\alpha^*\mc{O}_{\mb{G}(1,4)}(1) = c_1(M)$, the morphism $\alpha$ is given by a subsystem of $|T|$. In particular, $\alpha$ is a finite morphism. 

\begin{rem2.5}
The vector bundle $M$ is an example of a  Lazarsfeld-Mukai bundle. Lazarsfeld-Mukai bundles were studied in \cite{Laz86}, \cite{GL87}, \cite{Mu93} and were used in many constructions in algebraic geometry ever since. 
\end{rem2.5}

In what follows we shall use the notation of \cite{KL72} for Schubert calculus. In this notation, given a flag $A_0 \subsetneq A_1 \subset \bP^4$ with $a_i = \dim A_i$, the symbols $\Omega(a_0,a_1) = \Omega(A_0,A_1)$ denote the Schubert variety of lines $L$ in $\bP^4$ such that $\dim(L \cap A_i) \ge i$. The class of the Schubert variety in the cohomology ring $H^*(\mb{G}(1,4), \mb{Z})$ is denoted by the same symbol as the variety itself.

\begin{lem6} \label{lem6}
Every element of $|T|$ is irreducible and reduced. 
\end{lem6}
\begin{proof}
Recall the ongoing assumption that $(d,g)\in \ms{C}(V_5) \backslash\{(12,7),(13,8),(14,10)\}$. 
Suppose $T \sim D_1+D_2$ with $D_1,D_2$ non-trivial and write $D_1 \sim aT+bC$ for some $a,b \in \mb{Z}$.  Since $T$ is very ample, 
\begin{equation} \label{ineq1}
0 < D_1\cdot T < T^2=10.
\end{equation}
If $C$ is not a component of $D_1$ and $D_2$, then $D_i \cdot C \ge 0$, and therefore, 
\begin{equation} \label{ineq2} 
0 \le D_1\cdot C \le T\cdot C.
\end{equation} 
If $C$ is a component of either $D_1$ or $D_2$, without loss of generality, we shall assume that $D_1 = C$ and $D_2 = T-C$.  

Assume $(d,g) =(7,0)$, then the inequalities (\ref{ineq1}) and (\ref{ineq2}) become $0 < 10a+7b < 10$ and $0 \le 7a-2b \le 7$, respectively. Since either $a > 0 $ or $1-a>0$, we may check that the two inequalities may not hold simultaneously. Therefore, it remains to consider the case when $D_1=C$ and $D_2=T-C$. In this case we may check that $\deg D_2 = 3$ and $D_2^2=-6$. This implies that $D_2$ must have a line as a component. However, using $\mrm{Pic}(S) = \mb{Z}T \oplus \mb{Z}C$, we may check that $S$ does not contain lines. The remaining cases can be handled analogously and we omit the details. 
\end{proof}

\begin{rem3}
By \cite{Laz86}, $T$ is Brill-Noether general. In \cite{AM14} one may find the classification of Knutsen K3 surfaces all of whose hyperplane sections are irreducible and reduced. 
\end{rem3}

\begin{lem7} \label{lem7}
The composition of $\alpha\colon S \ra \mb{G}(1,4)$ with the Pl\"ucker embedding $\mb{G}(1,4) \hookrightarrow \bP^9$ has degree one.
\end{lem7}
\begin{proof} By the basis theorem \cite[p.1071]{KL72}, the class of $\alpha(S)$ in $H^*(\mb{G}(1,4), \mb{Z})$ can be expressed as
$[\alpha(S)]=a\Omega(0,3)+b\Omega(1,2)$
for some $a,b \in \mb{Z}$. Therefore, we may compute
$$4=\alpha_*\big(c_2(M)\big) = \alpha_*\alpha^*\Omega(2,3) = \deg(\alpha)\big(\Omega(2,3)\cdot [\alpha(S)]\big)=\deg(\alpha)b.$$
Furthermore, since $c_1(M)^2=T^2 = 10$, the degree of $\alpha$ must divide $10$. Hence, we have only two possibilities: (1) $\deg(\alpha)=1, a=6, b=4$; (2) $\deg(\alpha)=2, a=3, b=2$.

Assume that $\deg(\alpha) = 2$. Let $\langle \alpha(S) \rangle$ denote the linear span of $\alpha(S)$ in the Pl\"ucker embedding.   Since $h^0(\mc{O}_S(T))=7$ and $|T|$ is very ample, then $\deg(\alpha)=2$ implies $3 \le \dim \langle \alpha(S) \rangle  \le 5$. We shall eliminate the possibilities $\dim \langle S \rangle = 3, 4,$ or $5$ and conclude that $\alpha$ has degree one. 

Suppose $\dim \langle \alpha(S) \rangle = 3$. We may check that $\langle \alpha(S) \rangle$ is the intersection of all quadrics containing $\alpha(S)$ and since $\mb{G}(1,4) \subset \bP^9$ is also an intersection of quadrics then $\langle \alpha(S) \rangle \subset \mb{G}(1,4)$. Therefore, $\langle \alpha(S) \rangle$ is a maximal linear subvariety of $\mb{G}(1,4)$ and its class in $H^*(\mb{G}(1,4), \mb{Z})$ must be $\Omega(0,4)$. Since $\deg(\alpha)=2$, the surface $\alpha(S)$ has degree $5$ and is the intersection of $\langle \alpha(S) \rangle$ with a quintic hypersurface in $\bP^9$. Therefore, 
$$[\alpha(S)] = 5\Omega(2,4)\cdot \Omega(0,4) = 5\Omega(0,3),$$ 
which is impossible because $[\alpha(S)] = 3\Omega(0,3) + 2\Omega(1,2)$.

Suppose $\dim \langle \alpha(S) \rangle = 4$. Let $Y$ be the intersection of all quadrics containing $\alpha(S) \subset \bP^9$. Then $Y \subset \langle \alpha(S) \rangle \cap \mb{G}(1,4)$. Since $\mb{G}(1,4)$ does not contain any linear spaces of dimension 4, $\langle \alpha(S) \rangle$ is not in $\mb{G}(1,4)$. Hence, $Y \ne \langle \alpha(S) \rangle$. Also, since $\alpha(S)$ has degree 5 and is in $\bP^4$, $Y \ne \alpha(S)$. This shows that $Y$ is a quadric hypersurface in $\bP^4$.  Since the divisor $\alpha(S) \subset Y$ has degree $5$, $Y$ must be singular. Let $\alpha(S) \dashrightarrow Q$ be the map induced by the projection $Y \dashrightarrow Q$ from a singular point on $Y$. The variety $Q \subset \bP^3$ is either a smooth quadric or a quadric cone. In any case, a general hyperplane section of $Q$ is reducible and its pull-back to $S$ gives a reducible hyperplane section of $S$. This contradicts Lemma \ref{lem6}.

Suppose $\dim \langle \alpha(S) \rangle = 5$. In this case the morphism $\alpha \colon S \ra \mb{G}(1,4) \subset \bP^9$ can be factored into the closed embedding $S \hookrightarrow \bP^6$ given by $|T|$ followed by the projection $\pi\colon \bP^6 \dashrightarrow \bP^5 \subset \bP^9$ from some point $P \in \bP^6$. The inverse image $\Lambda:=\pi^{-1}(\alpha(S))$ is a cone whose vertex we shall denote by $P$. Since $S \subset \bP^6$ is the intersection of quadrics, there is a quadric $Q \subset \bP^6$ such that $S = \Lambda \cap Q$ (a priori $S \subset \Lambda \cap Q$, but $\deg S = \deg(\Lambda \cap Q) = 10$ implies the equality). Since $S$ is smooth, $\Lambda$ must be smooth along $S$. This implies that $\alpha(S)$ is smooth. Indeed, if $\alpha(S)$ were singular, $\Lambda$ would contain a line of singularities, which necessarily meets $Q$, giving a singular point of $\Lambda$ which lies on $S$. Therefore, $\alpha(S) \subset \bP^5$ is a smooth non-degenerate surface of degree $5$, hence  must be a del Pezzo surface. The Picard numbers of $S$ and $\alpha(S)$ are 2 and 5, respectively. This is impossible because $S \ra \alpha(S)$ is finite of degree 2.  
\end{proof}

\begin{lem5} \label{lem5}
The surface $\alpha(S)$ is not contained in any Schubert variety of type $\Omega(2,4)$ in $\mb{G}(1,4)$. 
\end{lem5}
\begin{proof} Suppose that $\alpha(S)$ is contained in the Schubert variety $\Omega(P_2,P_4)$, consisting of lines that meet a fixed $2$-plane $P_2 \subset P_4=\bP^4$. Consider the diagram 
$$\xymatrix{\bP(M) \ar[rd]^-\pi \ar[d]_-p& \\ S & P_4, }$$
where $p$ is the natural projection and $\pi$ is the morphism given by the complete linear system $|\mc{O}_{\bP(M)}(1)|$. Let $B$ be the image of $\pi$. 

Let $S^\p$ be the blow-up of $S$ along four general points that belong to a $g^1_4$ on $T\subset S$. We may check that $\pi(S^\p)$ is a hyperplane section of $B$.  Since $\alpha(S) \subset \Omega(P_2,P_4)$, the image of every fiber of $p$ under $\pi$ intersects $P_2$. There are two cases to consider. 

First, if $P_2 \subset B$ then $L:=P_2 \cap \pi(S^\p)$ is a line. If $\pi(S^\p)$ is general then $\pi^{-1}(L)$ is not a union of fibers of $p$ (otherwise, $B\subset P_2$, which is impossible, because $c_1(\mc{O}_{\mb{P}(M)}(1))^3=6$ and $\dim B=3$). Thus, we may and shall assume that  $\pi^{-1}(L)$ is not a union of fibers of $p$. A general hyperplane through $L$ cuts $\pi(S^\p)$ in a reducible curve, whose proper transform on $S$ is a reducible hyperplane section of $S$. This contradicts Lemma \ref{lem6}. 

Second, $P_2 \cap B$ is a plane curve $\Gamma$. If the images of all the fibers of $p$ pass through a single point, then $\alpha(S)$ is contained in a Schubert variety of type $\Omega(0,4)$, which is a $\bP^3$ in the Pl\"ucker embedding. By an argument as in the third paragraph of the proof of Lemma \ref{lem7} this is impossible because $[\alpha(S)] = 6\Omega(0,3) + 4\Omega(1,2)$. Thus, the images of the fibers of $p$ do not pass through a single point. Therefore, there is an irreducible component of $\Gamma$, which we also denote by $\Gamma$, such that for every $x \in \Gamma$, $p(\pi^{-1}(x))$ is a curve in $S$. A general hyperplane section of $B$ can be written as $\pi(S^\p)$, where $S^\p$ is as above. Let $x \in \Gamma \cap \pi(S^\p)$ and take a general hyperplane section $H$ of $\pi(S^\p)$ passing through $x$ (here general means that $H$ does not contain any of the images of the four exceptional divisors of the blow-up $S^\p \ra S$). The proper transform of $H$ in $S$ is a reducible hyperplane section of $S$, which contradicts Lemma \ref{lem6}.   
\end{proof}

\begin{pro5} \label{pro5}
The composition of $\alpha\colon S \ra \mb{G}(1,4)$ with the Pl\"ucker embedding $\mb{G}(1,4) \hookrightarrow \bP^9$ is given by the complete linear system $|T|$. In particular, the linear span $\langle \alpha(S) \rangle$ in $\bP^9$ is isomorphic to $\bP^6$. 
\end{pro5}
\begin{proof} By Lemma \ref{lem7}, the morphism $S \ra \bP^9$ has degree one. Therefore, it suffices to show that $\langle \alpha(S) \rangle$ has dimension 6. Suppose $\dim \langle \alpha(S) \rangle \le 5$, then $\alpha(S)$ is contained in a 3-dimensional family of hyperplanes in $\bP^9$. The subvariety of $(\bP^9)^*$  parametrizing hyperplanes $H \subset \bP^9$ such that $H \cap \mb{G}(1,4)$ is of type $\Omega(2,4)$ has dimension 6. By B\'ezout's theorem, this implies that $\alpha(S)$ is contained in a Schubert variety of type $\Omega(2,4)$, which contradicts Lemma \ref{lem5}.
\end{proof}

To simplify the notation in the sequel, we shall denote the image of $S$ in $\bP^9$ by $S$ as well. Also, the linear span of $S$ in $\bP^9$ will be denoted by $\langle S\rangle$.

\begin{cor5} \label{cor5}
The complete intersection of the linear span $\langle S \rangle$ with $\mb{G}(1,4)$ in $\bP^9$ is a smooth del Pezzo threefold $V_5$. 
\end{cor5}
\begin{proof} By Proposition \ref{pro5}, $\langle S \rangle \simeq \bP^6$. First, let us show that  
if $\mb{G}(1,4) \cap \langle S \rangle$ is singular, then there is a hyperplane $H \subset  \bP^9$ containing $\langle S \rangle$ such that $\mb{G}(1,4) \cap H$ is also singular. Let $V=\langle S \rangle \cap \mb{G}(1,4)$ and suppose $p$ is a singular point of $V$. Since $\dim V = 3$, the tangent space $T_pV$ has dimension $\ge 4$. Assume that $\dim T_pV = 4$. Let $v_1, v_2 \in T_p\mb{G}(1,4)$ be such that the linear span $\langle T_pV, v_1, v_2\rangle$ is all of $T_p\mb{G}(1,4)$. Let $l_1, l_2, l_3$ be the linear forms on $\bP^9$, whose zero locus is $\langle S \rangle$. The linear forms $l_i$ vanish on $T_pV$. Furthermore, we may find scalars $a,b,c$ such that the linear form $l:=al_1+bl_2+cl_3$ vanishes on $v_1$ and $v_2$, hence on all of $T_p\mb{G}(1,4)$. The zero locus $H:=Z(l)$ is the desired hyperplane. When $\dim T_pV >4$, the construction of $H$ is similar. 

Suppose $H$ is a hyperplane containing $\langle S \rangle$ such that $\mb{G}(1,4) \cap H$ is  singular. We shall show that $\mb{G}(1,4) \cap H$ is a Schubert variety of type $\Omega(2,4)$. The hyperplane section $H$ is given by a skew-symmetric $5\times 5$ matrix $A$. We may check that the singular locus of $\mb{G}(1,4) \cap H$ is the Grassmannian $\mb{G}(1, \bP(\ker A))$. Since $\mb{G}(1,4) \cap H$ is singular by assumption, $\dim \ker A = 3$. We may check that $\mb{G}(1,4) \cap H = \Omega(\bP(\ker A), \bP^4)$.

By Lemma \ref{lem5}, the image of $S$ is not contained in any hyperplane section of type $\Omega(2,4)$. Therefore, the intersection $\langle S \rangle \cap \mb{G}(1,4)$ is smooth and the lemma is proved. 
\end{proof}

\begin{thm1} \label{thm1}
For each $(d,g) \in \ms{C}(V_5) \backslash \{(14,10)\}$, the associated weak Fano threefold with small anti-canonical map exists. The numerical case with $(d,g) = (14,10)$ is not realizable. 
\end{thm1}

\begin{proof} Let us first show that the case with $(d,g)=(14,10)$ is not realizable. Suppose to the contrary that there is a weak Fano threefold $X$ obtained by blowing up a smooth curve $C$ of degree $14$ and genus $10$. A general member $\ti{S} \in |{-}K_X|$ is a smooth surface, whose image in $V_5$ is a smooth K3 surface $S$ containing $C$. By \cite[Prop.2.5, p.451]{JPR05}, $|{-}K_X|$ is free, and therefore, the linear system $|2T-C|$ is also free on $S$ (by the same reasoning as in the proof of Proposition \ref{pro3}). A general member $C^\p \in |2T-C|$ is a smooth curve of degree $6$ and genus $2$. 
The curve $C^\p$ lies on $V_5$ and the linear span of $C^\p$ in $\bP^6$ has dimension at most $4$. Thus, $C^\p$ is contained in two distinct hyperplanes $H_1$ and $H_2$ in $\bP^6$. Since $V_5$ is linearly normal and $\rho(V_5)=1$, the intersection $V_5\cap H_1\cap H_2$ is a curve of degree 5, which contains $C^\p$. However, $\deg C^\p=6$ and we have reached a contradiction. This shows that the case with $(d,g) = (14,10)$ is not realizable. 

Next, let us show that the cases with $(d,g) = (12,7)$ and $(13,8)$ are geometrically realizable. By \cite[A4, No.10]{JPR05} there is a weak Fano threefold obtained by blowing up a smooth curve $C^\p \subset V_5$ of degree $8$ and genus $3$. By the argument as in the previous paragraph, the curve $C^\p$ lies on a smooth K3 surface $S \subset V_5$ and the linear system $|2T-C^\p|$ is free. A general member $C \in |2T-C^\p|$ is a smooth curve of degree 12 and genus 7. The linear system $|2T-C|$ has no base points outside of its fixed components by \cite[Cor.3.2, p.611]{S-D74}. Also, since  $C^\p \in |2T-C|$ and $C^\p$ is not rational, it follows from \cite[2.7, p.610]{S-D74} that $|2T-C|$  is free. This implies that the blow-up of $V_5$ along $C$ is a weak Fano threefold, whose anti-canonical map must be small, because this numerical case does not appear on the tables in \cite{JPR05}. The case $(d,g)=(13,8)$ can be handled in a similar way by starting with the weak Fano threefold appearing as No.7 on Table E1-E5 in \cite{CM13}. 

In the remaining cases $(d,g)\in \ms{C}(V_5) \backslash\{(12,7),(13,8),(14,10)\}$. By Corollary \ref{cor5}, there is a smooth del Pezzo threefold $V_5$ containing Knutsen's K3 surface $S_{d,g}$ such that $C_{d,g} \subset S_{d,g}$. By Lemma \ref{lem9} the blow-up of $C_{d,g}$ on $V_5$ gives a weak Fano threefold with small anti-canonical map and $\phi^+$ is of type E1. 
\end{proof}

\begin{lem9} \label{lem9}
For each  $(d,g) \in \ms{C}(V_5) \backslash\{(12,7),(13,8),(14,10)\}$ the blow-up of $V_5$ at $C_{d,g}$ is a weak Fano threefold with small anti-canonical map and $\phi^+$ is of type E1. 
\end{lem9} 
\begin{proof} Throughout the proof we fix $(d,g)$ as in the statement of the lemma and omit the subscript $d,g$ from the notation $C_{d,g}$ and $S_{d,g}$ for simplicity. We let $X$ be the blow-up of $C$ on $V_5$ as before. Since $\mrm{Bs}(|{-}K_X|) \subset \mrm{Bs}(|2T-C|)$ then  by analogy with the proof of Proposition \ref{pro3}, to show that ${-}K_X$ is nef, it suffices to prove that $|2T-C|$ is free on $S$. Using Proposition \ref{pro2},  we may check that $|2T-C|$ is free on $S$ for the numerical cases listed in the lemma. 

It remains to check that $|{-}K_X|$ determines a small contraction and $\phi^+$ is of type E1. The numerical possibilities with $|{-}K_X|$ giving a divisorial contraction are classified in \cite[Table A.4, p.629]{JPR05}. If the morphism $\phi^+$ were not of type E1, then the listed numerical possibilities would have appeared either in \cite[7.4, 7.7, p.486]{JPR11} or in the non-E1-E1 tables in \cite{CM13}. Since this is not the case, $|{-}K_X|$ determines a small contraction and $\phi^+$ is of type E1.
\end{proof}

\section{Blow-ups of index one Fano threefolds} \label{X2n}

In this section we consider the cases on index one Fano threefolds $X_{2n}$ of genus $n+1$, where $n \in \{4,5,6,7,8,9,11\}$ (the customary notation for $X_{2n}$ is $X_{2g-2}$, where $g$ is the genus of $X_{2g-2}$, but in our case $g$ is reserved to denote the genus of the curve $C$). As in Section \ref{V5}, we may check that for each $(n,d,g)$ with $n$ as above and $(d,g) \in \mathscr{C}(X_{2n})$ there is a Knutsen K3 surface $S_{n,d,g} \subset \mb{P}^{n+1}$ with $\mrm{Pic}(S_{n,d,g}) = \mb{Z}T \oplus \mb{Z} C_{n,d,g}$, where $T$ is a general hyperplane section of $S_{n,d,g}$ (a smooth canonical curve of genus $n+1$) and $C_{n,d,g}$ is a smooth curve of degree $d$ and genus $g$. The existence of such surfaces was established in \cite{Knu02}. When there is no danger of confusion, we shall omit the subscript $n,d,g$ from $S_{n,d,g}$ and $C_{n,d,g}$ for simplicity. 

In the remainder of this section we shall use the following terminology introduced in \cite{Mu02}. By \cite[Def.3.8]{Mu02}, a polarized K3 surface $(S,T)$ is said to be \emph{Brill-Noether general} if $$h^0(M)h^0(N) < h^0(T)$$ for any non-trivial divisors $M,N$ such that $M + N \sim T$. It was shown in \cite{Mu02} that Brill-Noether general polarized K3 surfaces of genus 6,7,8,9,10 have projective models as complete intersections in certain homogeneous varieties, and as a consequence, are contained in a smooth Fano threefold of the same genus, Picard number 1 and index 1. Also, in \cite[Thm.5.5]{Mu02} Brill-Noether general K3 surfaces of genus $12$ are described as subvarieties of a Grassmannian in a natural way, and in particular, they are contained in a smooth Fano threefold of genus $12$, Picard number $1$ and index $1$. 

In the sequel we shall also use the result \cite[Thm.3.2]{Knu13} that classifies triples $(n,d,g)$ with $4 \le n \le 9$ for which there exists a Brill-Noether general K3 surface of degree $2n$ in $\bP^{n+1}$ containing a smooth irreducible curve of degree $d$ and genus $g$. 

\begin{pro6} \label{pro6}
For each $n \in \{4,5,6,7,8,9,11\}$ and $(d,g) \in \mathscr{C}(X_{2n})$ such that $(n,d,g)$ does not belong to the list
$$\{(6,4,1),(6,6,2),(7,7,2),(8,9,3),(9,8,2),(9,10,3),(9,11,4),(11,13,4),(11,14,5)\}$$
the blow-up of $X_{2n}$ at $C_{n,d,g}$ is a weak Fano threefold $X=X_{n,d,g}$. 
\end{pro6} 
\begin{proof}
Since $\mrm{Bs}(|{-}K_X|) \subset \mrm{Bs}(|T-C|)$, to show that ${-}K_X$ is nef, it suffices to prove that $|T-C|$ is free on $S$. Using Proposition \ref{pro2}, we may check that $|T-C|$ is free in the numerical cases listed in the proposition.  
\end{proof}

\begin{thm4} \label{thm4}
The numerical invariants listed in cases 3-8, 11-15, 18-20, 23, 55-57, 79, 93, 100, 104 in \cite[Table E1-E1]{CM13} of Sarkisov links starting with a Fano threefold of index 1 are geometrically realizable. The numerical cases 16, 17, 21, 22, 24-26, 58, 60 are not realizable. 
\end{thm4}
\begin{proof}
If $n \in \{4,5,6,7,8,9\}$ then using \cite[Thm.3.2]{Knu13} we may check that the Knutsen K3 surface $S_{n,d,g}$ (called \emph{Picard minimal $(n,d,g)$-surface} in \cite{Knu13}) with $(d,g)\in \mathscr{C}(X_{2n})$ is Brill-Noether general if and only if $(n,d,g)$ does not belong to the list 
$$\mathscr{N}=\{(6,4,1),(6,6,2),(7,7,2),(8,9,3),(9,8,2),(9,10,3),(9,11,4)\}.$$ Furthermore, by \cite[Thm.3.2]{Knu13} for $(n,d,g) \in \mathscr{N}$ there does not exist a Brill-Noether general K3 surface of degree $2n$ in $\bP^{n+1}$ containing a smooth irreducible curve of degree $d$ and genus $g$. This shows that the cases 16, 17, 21, 22, 24, 58, 60 are not geometrically realizable. 

For $(d,g) \in \mathscr{C}(X_{22})$ with $g=0,1$ we may check that the Knutsen K3 surface $S=S_{n,d,g}$ is Brill-Noether general by \cite[Lem.3.7]{Knu13}. For the cases 25, 26 with $(d,g)=(13,4),(14,5) \in \mathscr{C}(X_{22})$ there is no Brill-Noether general K3 surface of degree $22$ in $\bP^{12}$ containing a smooth curve of degree $d$ and genus $g$, because for the non-trivial decomposition $T=(T-C)+C$ we have $h^0(T-C)h^0(C) > h^0(T)$. For the remaining 3 cases $(d,g)=(10,2),(11,2),(12,3)$, we may check as in \cite[Lem.3.11]{Knu13} that $S_{n,d,g}$ is Brill-Noether general.

For each $(n,d,g)$ with $(d,g) \in \mathscr{C}(X_{2n})$ and such that $S_{n,d,g}$ is Brill-Noether general, there exists a smooth Fano threefold $X_{2n}$ containing $S_{n,d,g}$ by \cite{Mu02}. By Proposition \ref{pro6}, the blow-up of $X_{2n}$ at $C_{n,d,g}$ gives rise to a weak Fano threefold.

It remains to check that $|{-}K_X|$ determines a small contraction and $\phi^+$ is of type E1. The numerical possibilities with $|{-}K_X|$ giving a divisorial contraction are classified in \cite[Table A.4, p.629]{JPR05}. None of the cases listed as realizable in the statement of the theorem appear in \cite[Table A.4, p.629]{JPR05}. Therefore, $|{-}K_X|$ gives a small contraction. If the morphism $\phi^+$ were not of type E1, then the listed numerical possibilities would have appeared either in \cite[7.4, 7.7, p.486]{JPR11} or in the non-E1-E1 tables in \cite{CM13}. Since this is not the case for any of the numerical links listed as realizable in the statement of the theorem, $\phi^+$ is of type E1.
\end{proof}

\begin{rem4} \label{rem4}
Cases 2,10,78 also appear in the table \cite[A4, p.629]{JPR05} and cases 59,61,80 also appear in the tables \cite[7.4, 7.7, p.486]{JPR11}. It appears to be a delicate question whether the examples constructed in Proposition \ref{pro6} give geometric realizations of these numerical links. The authors are currently investigating these cases and hope to resolve them in a future paper.  
\end{rem4}

\begin{rem5} \label{rem5}
Case 28 of a numerical link starting with the del Pezzo threefold $V_2$ (which has index 2) appears to be unapproachable with the methods of this article due to the fact that $V_2$ does not admit a projective model as a complete intersection in a homogeneous space.
\end{rem5}

\section*{Acknowledgements} The authors express their gratitude to V.\,V.\,Shokurov for posing the problem and for his attention to this project. The authors are also grateful to Yu.\,G.\,Prokhorov, O.\,Debarre, S.\,I.\,Khashin, and A.\,Ortega for helpful correspondence and to the referee for suggesting Proposition \ref{pro2} that greatly simplified the exposition.

\section{Appendix}
The following tables summarize the numerical cases that were left open in \cite{CM13}. Case numbers are the same as those in \cite{CM13} but the cases are sorted by $-K_Y^3$ and the index of $Y$.

\label{E1E1table}
\begin{table}[h]
\caption{E1-E1}
\begin{center}
\begin{tabular}{|c|c|c|c|c|c|c|c|c|c|c|c|c|c|c|}
\hline
\textit{No.} & $-K_X^3$ & $-K_Y^3$ & $-K_{Y^+}^3$ & $\alpha$ & $\beta$ & $r$ & $d$ & $g$ & $r^+$ & $d^+$ & $g^+$ & $e/r^3$ & Exist? & Ref \\ \hline \hline
$\textit{2.}$ & 2 & 8 & 8 & 4 & -1 & 1 & 2 & 0 & 1 & 2 & 0 & 88 & ? & \ref{rem4}  \\ \hline


$\textit{3.}$ & 2 & 10 & 10 & 5 & -1 & 1 & 3 & 0 & 1 & 3 & 0 & 153 & :) & \ref{thm4} \\ \hline
$\textit{10.}$ & 2 & 10 & 10 & 4 & -1 & 1 & 4 & 1 & 1 & 4 & 1 & 56 & ? & \ref{rem4} \\ \hline

$\textit{4.}$ & 2 & 12 & 12 & 6 & -1 & 1 & 4 & 0 & 1 & 4 & 0 & 248 & :) & \ref{thm4}  \\ \hline
$\textit{11.}$ & 2 & 12 & 12 & 5 & -1 & 1 & 5 & 1 & 1 & 5 & 1 & 115 & :) & \ref{thm4} \\ \hline
$\textit{16.}$ & 2 & 12 & 12 & 4 & -1 & 1 & 6 & 2 & 1 & 6 & 2 & 24 & x & \ref{thm4}  \\ \hline
$\textit{58.}$ & 4 & 12 & 12 & 2 & -1 & 1 & 4 & 1 & 1 & 4 & 1 & 8 & x & \ref{thm4} \\ \hline
$\textit{5.}$ & 2 & 14 & 14 & 7 & -1 & 1 & 5 & 0 & 1 & 5 & 0 & 379 & :) & \ref{thm4} \\ \hline
$\textit{12.}$ & 2 & 14 & 14 & 6 & -1 & 1 & 6 & 1 & 1 & 6 & 1 & 204 & :) & \ref{thm4} \\ \hline
$\textit{17.}$ & 2 & 14 & 14 & 5 & -1 & 1 & 7 & 2 & 1 & 7 & 2 & 77 & x & \ref{thm4} \\ \hline
$\textit{55.}$ & 4 & 14 & 14 & 3 & -1 & 1 & 4 & 0 & 1 & 4 & 0 & 68 & :) & \ref{thm4} \\ \hline
$\textit{6.}$ & 2 & 16 & 16 & 8 & -1 & 1 & 6 & 0 & 1 & 6 & 0 & 552 & :) &  \ref{thm4}\\ \hline
$\textit{13.}$ & 2 & 16 & 16 & 7 & -1 & 1 & 7 & 1 & 1 & 7 & 1 & 329 & :) & \ref{thm4}  \\ \hline
$\textit{18.}$ & 2 & 16 & 16 & 6 & -1 & 1 & 8 & 2 & 1 & 8 & 2 & 160 & :) & \ref{thm4} \\ \hline
$\textit{21.}$ & 2 & 16 & 16 & 5 & -1 & 1 & 9 & 3 & 1 & 9 & 3 & 39 & x & \ref{thm4}  \\ \hline
$\textit{59.}$ & 4 & 16 & 16 & 3 & -1 & 1 & 6 & 1 & 1 & 6 & 1 & 42 & ? &  \ref{rem4}\\ \hline
$\textit{78.}$ & 6 & 16 & 16 & 2 & -1 & 1 & 4 & 0 & 1 & 4 & 0 & 32 & ? & \ref{rem4}  \\ \hline

$\textit{7.}$ & 2 & 18 & 18 & 9 & -1 & 1 & 7 & 0 & 1 & 7 & 0 & 773 & :) &  \ref{thm4} \\ \hline
$\textit{14.}$ & 2 & 18 & 18 & 8 & -1 & 1 & 8 & 1 & 1 & 8 & 1 & 496 & :) & \ref{thm4}  \\ \hline
$\textit{19.}$ & 2 & 18 & 18 & 7 & -1 & 1 & 9 & 2 & 1 & 9 & 2 & 279 & :) & \ref{thm4} \\ \hline
$\textit{22.}$ & 2 & 18 & 18 & 6 & -1 & 1 & 10 & 3 & 1 & 10 & 3 & 116 & x & \ref{thm4}  \\ \hline
$\textit{24.}$ & 2 & 18 & 18 & 5 & -1 & 1 & 11 & 4 & 1 & 11 & 4 & 1 & x &  \ref{thm4}\\ \hline
$\textit{56.}$ & 4 & 18 & 18 & 4 & -1 & 1 & 6 & 0 & 1 & 6 & 0 & 144 & :) & \ref{thm4} \\ \hline
$\textit{60.}$ & 4 & 18 & 18 & 3 & -1 & 1 & 8 & 2 & 1 & 8 & 2 & 16 & x & \ref{thm4} \\ \hline
$\textit{80.}$ & 6 & 18 & 18 & 2 & -1 & 1 & 6 & 1 & 1 & 6 & 1 & 12 & ? & \ref{rem4} \\ \hline
$\textit{100.}$ & 10 & 18 & 18 & 1 & -1 & 1 & 3 & 0 & 1 & 3 & 0 & 9 & :) & \ref{thm4}  \\ \hline

\end{tabular}
\label{tb:E1E1Table1}
\end{center}
\end{table}

\begin{table}
\caption{E1-E1 (continued)}
\begin{center}
\begin{tabular}{|c|c|c|c|c|c|c|c|c|c|c|c|c|c|c|}
\hline
\textit{No.} & $-K_X^3$ & $-K_Y^3$ & $-K_{Y^+}^3$ & $\alpha$ & $\beta$ & $r$ & $d$ & $g$ & $r^+$ & $d^+$ & $g^+$ & $e/r^3$ & Exist? & Ref \\ \hline \hline

$\textit{8.}$ & 2 & 22 & 22 & 11 & -1 & 1 & 9 & 0 & 1 & 9 & 0 & 1383 & :) &  \ref{thm4} \\ \hline
$\textit{15.}$ & 2 & 22 & 22 & 10 & -1 & 1 & 10 & 1 & 1 & 10 & 1 & 980 & :) & \ref{thm4} \\ \hline
$\textit{20.}$ & 2 & 22 & 22 & 9 & -1 & 1 & 11 & 2 & 1 & 11 & 2 & 649 & :) & \ref{thm4} \\ \hline
$\textit{23.}$ & 2 & 22 & 22 & 8 & -1 & 1 & 12 & 3 & 1 & 12 & 3 & 384 & :) & \ref{thm4} \\ \hline
$\textit{25.}$ & 2 & 22 & 22 & 7 & -1 & 1 & 13 & 4 & 1 & 13 & 4 & 179 & x & \ref{thm4} \\ \hline
$\textit{26.}$ & 2 & 22 & 22 & 6 & -1 & 1 & 14 & 5 & 1 & 14 & 5 & 28 & x & \ref{thm4} \\ \hline
$\textit{57.}$ & 4 & 22 & 22 & 5 & -1 & 1 & 8 & 0 & 1 & 8 & 0 & 268 & :) & \ref{thm4}  \\ \hline
$\textit{61.}$ & 4 & 22 & 22 & 4 & -1 & 1 & 10 & 2 & 1 & 10 & 2 & 80 & ? & \ref{rem4} \\ \hline
$\textit{79.}$ & 6 & 22 & 22 & 3 & -1 & 1 & 7 & 0 & 1 & 7 & 0 & 89 & :) & \ref{thm4}  \\ \hline
$\textit{93.}$ & 8 & 22 & 22 & 2 & -1 & 1 & 6 & 0 & 1 & 6 & 0 & 36 & :) & \ref{thm4}  \\ \hline
$\textit{104.}$ & 12 & 22 & 22 & 1 & -1 & 1 & 4 & 0 & 1 & 4 & 0 & 8 & :) & \ref{thm4} \\ \hline

$\textit{28.}$ & 2 & 16 & 16 & 8 & -1 & 2 & 3 & 0 & 2 & 3 & 0 & 69 & ? & \ref{rem5}\\ \hline

$\textit{30.}$ & 2 & 32 & 32 & 16 & -1 & 2 & 7 & 0 & 2 & 7 & 0 & 521 & :) & \ref{thm3} \\ \hline
$\textit{34.}$ & 2 & 32 & 32 & 14 & -1 & 2 & 8 & 2 & 2 & 8 & 2 & 328 & :) & \ref{thm3} \\ \hline
$\textit{37.}$ & 2 & 32 & 32 & 12 & -1 & 2 & 9 & 4 & 2 & 9 & 4 & 183 & :) & \ref{thm3}  \\ \hline
$\textit{39.}$ & 2 & 32 & 32 & 10 & -1 & 2 & 10 & 6 & 2 & 10 & 6 & 80 & ? & \ref{rem2} \\ \hline
$\textit{41.}$ & 2 & 32 & 32 & 8 & -1 & 2 & 11 & 8 & 2 & 11 & 8 & 13 & :) & \ref{thm3}  \\ \hline
$\textit{64.}$ & 4 & 32 & 18 & 3.5 & -0.5 & 2 & 7 & 1 & 1 & 7 & 1 & 77 & :) & \ref{thm3} \\ \hline
$\textit{67.}$ & 4 & 32 & 32 & 6 & -1 & 2 & 8 & 3 & 2 & 8 & 3 & 40 & ? &   \ref{rem2} \\ \hline
$\textit{84.}$ & 6 & 32 & 32 & 4 & -1 & 2 & 7 & 2 & 2 & 7 & 2 & 17 & :) & \ref{thm3} \\ \hline

$\textit{31.}$ & 2 & 40 & 40 & 20 & -1 & 2 & 9 & 0 & 2 & 9 & 0 & 1011 & :) &  \ref{thm1}  \\ \hline
$\textit{35.}$ & 2 & 40 & 40 & 18 & -1 & 2 & 10 & 2 & 2 & 10 & 2 & 710 & :) &  \ref{thm1} \\ \hline
$\textit{38.}$ & 2 & 40 & 40 & 16 & -1 & 2 & 11 & 4 & 2 & 11 & 4 & 469 & :) &  \ref{thm1}  \\ \hline
$\textit{40.}$ & 2 & 40 & 40 & 14 & -1 & 2 & 12 & 6 & 2 & 12 & 6 & 282 & :) &  \ref{thm1}  \\ \hline
$\textit{42.}$ & 2 & 40 & 40 & 12 & -1 & 2 & 13 & 8 & 2 & 13 & 8 & 143 & :) &  \ref{thm1}  \\ \hline
$\textit{43.}$ & 2 & 40 & 40 & 10 & -1 & 2 & 14 & 10 & 2 & 14 & 10 & 46 & x &  \ref{thm1} \\ \hline
$\textit{65.}$ & 4 & 40 & 22 & 4.5 & -0.5 & 2 & 9 & 1 & 1 & 9 & 1 & 171 & :) &  \ref{thm1} \\ \hline
$\textit{68.}$ & 4 & 40 & 40 & 8 & -1 & 2 & 10 & 3 & 2 & 10 & 3 & 110 & :) &  \ref{thm1} \\ \hline
$\textit{69.}$ & 4 & 40 & 40 & 6 & -1 & 2 & 12 & 7 & 2 & 12 & 7 & 18 & :) &  \ref{thm1} \\ \hline
$\textit{81.}$ & 6 & 40 & 18 & 2.5 & -0.5 & 2 & 9 & 2 & 1 & 5 & 0 & 47 & :) &  \ref{thm1} \\ \hline
$\textit{83.}$ & 6 & 40 & 40 & 6 & -1 & 2 & 8 & 0 & 2 & 8 & 0 & 82 & :) &  \ref{thm1}  \\ \hline

\end{tabular}
\label{tb:E1E1Table2}
\end{center}
\end{table}

\begin{table}
\caption{E1-E1 (continued)}
\begin{center}
\begin{tabular}{|c|c|c|c|c|c|c|c|c|c|c|c|c|c|c|}
\hline
\textit{No.} & $-K_X^3$ & $-K_Y^3$ & $-K_{Y^+}^3$ & $\alpha$ & $\beta$ & $r$ & $d$ & $g$ & $r^+$ & $d^+$ & $g^+$ & $e/r^3$ & Exist? & Ref \\ \hline \hline

$\textit{94.}$ & 8 & 40 & 16 & 1.5 & -0.5 & 2 & 9 & 3 & 1 & 3 & 0 & 12 & :) &  \ref{thm1}  \\ \hline
$\textit{96.}$ & 8 & 40 & 40 & 4 & -1 & 2 & 8 & 1 & 2 & 8 & 1 & 28 & :) &  \ref{thm1}  \\ \hline
$\textit{101.}$ & 10 & 40 & 22 & 1.5 & -0.5 & 2 & 7 & 0 & 1 & 5 & 0 & 18 & :) &  \ref{thm1}  \\ \hline

$\textit{44.}$ & 2 & 54 & 54 & 25 & -1 & 3 & 9 & 2 & 3 & 9 & 2 & 571 & :) & \ref{thm2} \\ \hline
$\textit{45.}$ & 2 & 54 & 54 & 22 & -1 & 3 & 10 & 5 & 3 & 10 & 5 & 372 & :) & \ref{thm2} \\ \hline
$\textit{46.}$ & 2 & 54 & 54 & 19 & -1 & 3 & 11 & 8 & 3 & 11 & 8 & 221 & :) & \ref{thm2} \\ \hline
$\textit{47.}$ & 2 & 54 & 54 & 16 & -1 & 3 & 12 & 11 & 3 & 12 & 11 & 112 & ? & \ref{rem1}\\ \hline
$\textit{48.}$ & 2 & 54 & 54 & 13 & -1 & 3 & 13 & 14 & 3 & 13 & 14 & 39 & :) & \ref{thm2} \\ \hline
$\textit{70.}$ & 4 & 54 & 16 & 11/3 & -1/3 & 3 & 9 & 3 & 1 & 5 & 0 & 103 & :) & \ref{thm2} \\ \hline
$\textit{71.}$ & 4 & 54 & 54 & 13 & -1 & 3 & 8 & 0 & 3 & 8 & 0 & 164 & :) & \ref{thm2} \\ \hline
$\textit{72.}$ & 4 & 54 & 54 & 10 & -1 & 3 & 10 & 6 & 3 & 10 & 6 & 60 & ? &  \ref{rem1} \\ \hline
$\textit{86.}$ & 6 & 54 & 22 & 8/3 & -1/3 & 3 & 8 & 1 & 1 & 8 & 1 & 48 & :) &\ref{thm2} \\ \hline
$\textit{88.}$ & 6 & 54 & 54 & 7 & -1 & 3 & 9 & 4 & 3 & 9 & 4 & 31 & :) & \ref{thm2} \\ \hline
$\textit{97.}$ & 8 & 54 & 18 & 5/3 & -1/3 & 3 & 8 & 2 & 1 & 4 & 0 & 20 & :) & \ref{thm2}  \\ \hline
$\textit{102.}$ & 10 & 54 & 54 & 4 & -1 & 3 & 8 & 3 & 3 & 8 & 3 & 8 & ? &   \ref{rem1} \\ \hline
$\textit{105.}$ & 12 & 54 & 40 & 7/3 & -2/3 & 3 & 7 & 1 & 2 & 7 & 1 & 7 & :) & \ref{thm2}  \\ \hline

\end{tabular}
\label{tb:E1E1Table3}
\end{center}
\end{table}


\begin{thebibliography}{99}
\bibitem[AM14]{AM14} Arap, M., Marshburn, N.: \emph{Brill-Noether general curves on Knutsen K3 surfaces.} Comptes rendus - Math\'ematique \textbf{352} (2014), 133--135.
\bibitem[Bea96]{Bea96} Beauville, A.: \emph{Complex algebraic surfaces.} Second edition. London Mathematical Society Student Texts, \textbf{34}. Cambridge University Press, 1996.
\bibitem[BL12]{BL12} Blanc, J., Lamy, S.: \emph{Weak Fano threefolds obtained by blowing-up a space curve and construction of Sarkisov links}. Proc. Lond. Math. Soc. (3) \textbf{105} (2012), no. 5, 1047--1075. 
\bibitem[CHNP12]{CHNP12}  Corti, A., Haskins, M., Nordstr\"om, J., Pacini, T.: \emph{$G_2$-manifolds and associative submanifolds via semi-Fano 3-folds.}  arXiv 1207.4470v3.
\bibitem[CHNP13]{CHNP13} Corti, A., Haskins, M., Nordstr\"om, J., Pacini, T.: \emph{Asymptotically cylindrical Calabi-Yau 3-folds from weak Fano 3-folds.} Geom. Topol. \textbf{17} (2013), 1955--2059.
\bibitem[CM13]{CM13} Cutrone, J.W., Marshburn, N.A.: \emph{Towards the classification of weak Fano threefolds with $\rho = 2$.} Cent. Eur. J. Math. \textbf{11} (2013), no. 9, 1552--1576.
\bibitem[GL87]{GL87} Green, M., Lazarsfeld, R.: \emph{Special divisors on curves on a  K3 surface}. Invent. math., Volume \textbf{89}(1987), Issue 2, 357--370.
\bibitem[Gu82]{Gu82} Gushel, N. P.: \emph{Fano varieties of genus $6$.} (Russian) Izv. Akad. Nauk SSSR Ser. Mat. \textbf{46} (1982), no. 6, 1159--1174,
\bibitem[Gu92]{Gu92} Gushel, N. P.: \emph{Fano $3$-folds of genus $8$.} (Russian) Algebra i
Analiz \textbf{4} (1992), no. 1, 120--134; translation in St. Petersburg Math. J. \textbf{4} (1993), no. 1, 115--129
\bibitem[Isk89]{Isk89} Iskovskikh, V. A.: \emph{Double projection from a line onto Fano $3$-folds of the first kind.} (Russian) Mat. Sb. \textbf{180} (1989), no. 2, 260--278; translation in Math. USSR-Sb. \textbf{66} (1990), no. 1, 265--284.
\bibitem[Isk77]{Isk77}Iskovskih, V. A.: \emph{Fano threefolds. I.} (Russian) Izv. Akad. Nauk SSSR Ser. Mat. \textbf{41}, no. 3 (1977),  516--562. 
\bibitem[Isk78]{Isk78}Iskovskih, V. A.: \emph{Fano threefolds. II.} (Russian) Izv. Akad. Nauk SSSR Ser. Mat. \textbf{42}, no. 3 (1978), 506--549. 
\bibitem[IP99]{IP99} Iskovskikh, V. A., Prokhorov, Yu. G.: \emph{Fano varieties. Algebraic
geometry V}, Encyclopaedia Math. Sci., \textbf{47}, Springer, Berlin, 1999.
\bibitem[JPR05]{JPR05} Jahnke, P., Peternell, T., Radloff, I.: \emph{Threefolds with big
and nef anti-canonical bundles I.} Math. Ann. \textbf{333} (2005), no. 3, 569--631.
\bibitem[JPR11]{JPR11} Jahnke, P., Peternell, T., Radloff, I.: \emph{Threefolds with big
and nef anti-canonical bundles II.} Cent. Eur. J. Math. \textbf{9} (2011), no. 3, 449--488.
\bibitem[KL72]{KL72} Kleiman, S. L., Laksov, D.: \emph{Schubert calculus.} Amer. Math. Monthly \textbf{79} (1972), 1061--1082.
\bibitem[Knu02]{Knu02} Knutsen, A.L.: \emph{Smooth curves on projective $K3$ surfaces.} Math. Scand. \textbf{90} (2002) 215--231.
\bibitem[Knu13]{Knu13} A.L. Knutsen: \emph{Smooth, isolated curves in families of Calabi-Yau threefolds in homogeneous spaces.} J. Korean Math. Soc. \textbf{50} (2013), No. 5, pp. 1033--1050.
\bibitem[Ko89]{Ko89} Koll\'ar, J.: \emph{Flops.} Nagoya Math. J. \textbf{113} (1989), 15--36.
\bibitem[LazI]{LazI} Lazarsfeld, R.: \emph{Positivity in algebraic geometry I. Classical
setting: line bundles and linear series.} Ergeb. der Math. und ihrer Grenz. 3. Folge, \textbf{48}. Springer-Verlag, Berlin, 2004.
\bibitem[Laz86]{Laz86} Lazarsfeld, R.: \emph{Brill-Noether-Petri without degenerations.} J.
Differential Geom. \textbf{23} (1986), no. 3, 299--307.
\bibitem[LeB82]{LeB82} Le Barz, P.: \emph{Formules multis\'ecantes pour les courbes gauches
quelconques.} (Nice, 1981), pp. 165--197, Progr. Math., \textbf{24}, BirkhŠuser, Boston, Mass., 1982.
\bibitem[Ma73]{Ma73} Maruyama, M.: \emph{On a family of algebraic vector bundles.} Number
theory, algebraic geometry and commutative algebra, in honor of Yasuo Akizuki, pp. 95Ð146. Kinokuniya, Tokyo, 1973.
\bibitem[Mo82]{Mo82} Mori, S.: \emph{Threefolds whose canonical bundles are not numerically
effective.} Ann. of Math. (2) \textbf{116} (1982), no. 1, 133--176.
\bibitem[MM84]{MM84} Mori, S., Mukai, S.: \emph{Classification of Fano 3-folds with $B_2 \ge 2$. I.} Algebraic and topological theories (Kinosaki, 1984), 496--545, Kinokuniya, Tokyo, 1986.
\bibitem[Mu93]{Mu93} Mukai, S.: \emph{Curves and Grassmannians.} Algebraic geometry and related topics (Inchon, 1992), Int. Press, Cambridge, MA (1993) 19--40.
\bibitem[Mu02]{Mu02}  Mukai, S.: \emph{New development of theory of Fano 3-folds: vector bundle method and moduli problem.} Sugaku \textbf{47} (1995), no. 2, 125--144; translation in: Sugaku Expositions \textbf{15} (2002), no. 2, 125--150.
\bibitem[Mu10]{Mu10} Mukai, S.: \emph{Curves and symmetric spaces, II.} Annals of Math., \textbf{172} (2010), 1539--1558.
\bibitem[S-D74]{S-D74} Saint-Donat, B.: \emph{Projective models of K3 surfaces.} Amer. J. Math. \textbf{96} (1974), 602--639.
\bibitem[Sh79]{Sh79} Shokurov, V. V.: \emph{The existence of a line on Fano varieties.} (Russian) Izv. Akad. Nauk SSSR Ser. Mat. \textbf{43} , no. 4 (1979), 922--964
\bibitem[Tak09]{Tak09} Takeuchi, K.:  \emph{Weak Fano threefolds with del Pezzo fibration.} arXiv preprint (2009).
\end{thebibliography}
\end{document}